\documentclass[reqno,12pt]{amsart}

\newtheorem{theorem}{Theorem}[section]
\newtheorem{lemma}[theorem]{Lemma}
\newtheorem{corollary}[theorem]{Corollary}

\theoremstyle{definition}
\newtheorem{assumption}[theorem]{Assumption}

\theoremstyle{remark}
\newtheorem{remark}[theorem]{Remark}

\makeatletter
\def\dashint{\operatorname%
{\,\,\text{\bf--}\kern-.98em\DOTSI\intop\ilimits@\!\!}}
\makeatother

\newcommand\bL{\mathbb{L}}
\newcommand\bM{\mathbb{M}}

\newcommand\bR{\mathbb{R}}

\newcommand\cF{\mathcal{F}}

\newcommand\cL{\mathcal{L}}
\newcommand\cM{\mathcal{M}}
\newcommand\cS{\mathcal{S}}

\newcommand\frB{\mathfrak{B}}
\newcommand\frL{\mathfrak{L}}
\newcommand\frU{\mathfrak{U}}

 \newcommand{\mysection}[1]{\section{#1}
 \setcounter{equation}{0}}
 
\def\bC{\mathbb{C}}
\def\bZ{\mathbb{Z}}

\begin{document}

\title[Bellman's equations with VMO coefficients]
{On Bellman's equations with VMO coefficients}
\author{N.V. Krylov}
\thanks{The work  was partially supported by
NSF Grant DMS-0653121}
\email{krylov@math.umn.edu}
\address{127 Vincent Hall, University of Minnesota,
 Minneapolis, MN, 55455}
 
 

\begin{abstract}
We present a result about solvability in
$W^{2}_{p}$, $p>d$, in the whole space $\bR^{d}$
of Bellman's equations with VMO ``coefficients''.
Parabolic equations are touched upon as well.
\end{abstract}

\maketitle

\mysection{Main result}

Let $\bR^{d}$ be the Euclidean space of 
points $x=(x^{1},...,x^{d})$, $x^{i}\in\bR=(-\infty,\infty)$.
Fix a $\delta\in(0,1)$ and denote by $\cS_{\delta}$ the set of
symmetric $d\times d$-matrices $a=(a^{ij})$ satisfying
$$
\delta|\xi|^{2}\leq a^{ij}\xi^{i}\xi^{j}\leq\delta^{-1}|\xi|^{2},
\quad\forall \xi\in\bR^{d}.
$$
Let $\Omega$ be a separable metric space
 and assume that for any
$\omega\in \Omega$ and $x\in\bR^{d}$ we are given
$a(\omega,x)\in\cS_{\delta}$,
$b(\omega,x)\in\bR^{d}$, and $c(\omega,x),
f(\omega,x)\in\bR$. We assume that 
these functions are measurable in $x$ for each $\omega$,
continuous in $\omega$ for each $x$, and 
$$
|b(\omega,x)|+c(\omega,x)\leq K,\quad c(\omega,x)\geq0,
\quad\forall \omega,x,
$$
$$
\bar{f}(x):=\sup_{\omega\in \Omega}|f(\omega,x)|<\infty
\quad\forall x,
$$
where $K$ is a fixed constant. Observe that, owing to the continuity
of $f$ in $\omega$ and separability of $\Omega$, the function $\bar{f}$
is measurable.
For $r>0$ and $x\in\bR^{d}$ set
$$
B_{r}(x)=\{y\in\bR^{d}:|x-y|<r\},\quad B_{r}=B_{r}(0).
$$
For a measurable set $\Gamma\subset\bR^{d}$ by
$|\Gamma|$ we denote its volume. In Section \ref{section 12.23.1}
we will use the same notation for measurable $\Gamma\subset\bR^{d+1}$.
Introduce,
$$
(u)_{\Gamma}=
\dashint_{\Gamma}u(x)\,dx=\frac{1}{|\Gamma|}\int_{\Gamma}
u(x)\,dx.
$$

In particular,
$$
(a)_{B_{r}(x)}(\omega )=\dashint_{B_{r}(x)}a(\omega,y)\,dy
$$

In the following assumption there is a parameter $\theta
\in(0,1]$, whose value will be specified later.
\begin{assumption}
                                        \label{assumption 12.18.1}
There exists an $R_{0}\in(0,\infty)$ such that
for any $r\in(0,R_{0}]$  and $x\in\bR^{d}$ we have
\begin{equation}
                                                    \label{12.24.6}
\dashint_{B_{r}(x)}\sup_{\omega\in \Omega} 
|a(\omega,y)-(a)_{B_{r}(x)}(\omega )|\,dy \leq\theta.
\end{equation}

\end{assumption}
Observe that if $a$ is independent of $\omega$ (semilinear equations),
and for any $\theta>0$ there is an $R_{0}>0$ such that \eqref{12.24.6}
is satisfied for any $r\in(0,R_{0}]$  and $x\in\bR^{d}$, then
$a\in VMO$.

For constant $\lambda>0$ we will be considering the following equation
$$
\sup_{\omega\in \Omega}[a^{ij}(\omega,x)D_{ij}u(x)
+b^{i}(\omega,x)D_{i}u(x)
$$
\begin{equation}
                                               \label{12.23.5}
-(c(\omega,x)+\lambda)u(x)+f(\omega,x)]=0
\end{equation}
in $\bR^{d}$. Of course, $D_{i}=\partial/\partial x^{i}$,
$D_{ij}=D_{i}D_{j}$. By $Du$ we denote the gradient of $u$
and $D^{2}u$ the Hessian matrix of $u$.
If $D$ is an open set and $p\geq1$, by $\cL_{p}(D)$ we denote the usual
Lebesgue space and by $W^{2}_{p}(D)$ the usual Sobolev space.
Denote $\cL_{p}=\cL_{p}(\bR^{d})$, $W^{2}_{p}=W^{2}_{p}
(\bR^{d})$.
 We say that a function $u\in W^{2}_{p}$
satisfies \eqref{12.23.5} in $\bR^{d}$ if \eqref{12.23.5} holds
almost everywhere in $\bR^{d}$.

The main information about the solvability
of \eqref{12.23.5} will be obtained while studying
its reduced form
\begin{equation}
                                               \label{12.6.1}
\sup_{\omega\in \Omega}[a^{ij}(\omega,x)D_{ij}u(x)+f(\omega,x)]=0.
\end{equation}
 
\begin{theorem}
                                   \label{theorem 12.23.3}
Let $p>d$. Then there exists a constant
$\theta\in(0,1)$ depending only on $p,d$, and $\delta$
such that, if Assumption \ref{assumption 12.18.1} holds with this $\theta$,
then there exists a $\lambda_{0}=\lambda_{0}
(\delta,p,K,d)\geq0$ such that for any $\lambda
\geq\lambda_{0}$ and
$u\in W^{2}_{p}$ satisfying \eqref{12.23.5}
we have
$$
\lambda\|u\|_{\cL_{p}}+\|D^{2}u\|_{\cL_{p}}
\leq N\|\bar{f}\|_{\cL_{p}},
$$
where $N=N(\delta,p, d)$.
Moreover,  for any $\lambda
>\lambda_{0}$ there is a unique solution of \eqref{12.23.5}
in $W^{2}_{p}$.

\end{theorem}

Previously the $W^{2}_{p}$-estimates for a class of fully nonlinear
elliptic equations were obtained by Caffarelli
in \cite{Ca} (see also \cite{CC}). It seems to the author that
the class of equations from \cite{CC}
 does not include Bellman's equations like \eqref{12.6.1} when
$\Omega$ consists of only two points:
$$
\max\big(a^{ij}_{1}(x)D_{ij}u(x) +f_{1}(x),
a^{ij}_{2}(x)D_{ij}u(x) +f_{2}(x)\big)=0
$$
unless $a_{1}$ and $a_{2}$ are uniformly sufficiently close
to  continuous functions and $f_{1}-f_{2}$ is uniformly sufficiently 
close to a uniformly continuous function.

Our methods are slightly different from those from \cite{CC}.
As in \cite{CC} we start from the result of \cite{FH} but then
we follow more closely the general scheme  from \cite{Kr08}.
We are only dealing with equations in the whole space
or interior estimates. One can probably obtain similar results
for equations in half spaces and smooth domains by following the
recent approach developed in \cite{DK} for {\em linear\/} equations
and systems. In \cite{DK} one can also find an extensive list of
references on linear equations with VMO coefficients.

We prove Theorem \ref{theorem 12.23.3} in Section
\ref{section 12.23.4} after we prepare some necessary
tools in the next section. In Section \ref{section 12.23.1}
we give some comments on how one can start
obtaining similar results for parabolic equations.
In the final Section \ref{section 12.24.3} we prove
the necessary facts from Real Analysis.

 The author is sincerely grateful to
Hongjie Dong for pointing out several
glitches in the first draft of the paper.

\mysection{Equations with constant coefficients}
                                                   \label{section 12.23.3}

In this section we consider equation \eqref{12.6.1}
under the additional assumption that $a^{ij} $
are independent of $x$, namely,
\begin{equation}
                                               \label{12.16.3}
\sup_{\omega\in \Omega}\big[
a^{ij} (\omega )D_{ij}u(x)+f(\omega,x)\big]=0.
\end{equation}

\begin{lemma}
                                         \label{lemma 12.6.1}
Let $u\in W^{2}_{d,loc}$. Then for any $r\in(0,\infty)$
$$
\sup_{B_{r}}|u(x)-x^{i}(u_{x^{i}})_{B_{r}}-(u )_{B_{r}}|^{d}\leq
N(d)r^{2d}\dashint_{B_{r}}|D^{2}u|^{d}\,dx.
$$

\end{lemma}

Proof. First assume that $r=1$
and set
$$
v(x)=u(x)-x^{i}(u_{x^{i}})_{B_{1}}-(u )_{B_{1}}
$$ 
 By Sobolev embedding theorems the
left-hand side is less than a constant times
$$
\sup_{B_{1}}|v|\leq N\|v\|_{W^{2}_{d}(B_{1})}
$$
and by Poincar\'e's inequality
$$
\|v\|_{W^{2}_{d}(B_{1})}^{d}=
\|D^{2}u\|_{\cL_{d}(B_{1})}^{d}+
\|D u-(Du)_{B_{1}}\|_{\cL_{d}(B_{1})}^{d}
$$
$$
+\|v\|_{\cL_{d}(B_{1})}^{d}\leq
N\|D^{2}u\|_{\cL_{d}(B_{1})}^{d}.
$$
For general $r$ it suffices to use dilations.
The lemma is proved.

\begin{lemma}
                                            \label{lemma 12.16.01}

Let $r\in(0,\infty)$, $\nu\geq2$ and let
 $v\in C^{2}(\bar{B}_{\nu r})$ be a solution of
\eqref{12.16.3} in $B_{\nu r}$ with $f \equiv0$.
Then there are  constants $\beta\in(0,1)$ and $N$,
depending only on $d$ and $\delta$, such that
$$
\dashint_{B_{r}}\dashint_{B_{r}}|D^{2}v(x)-D^{2}v(y)|
\,dxdy\leq N\nu^{-2-\beta} r^{-2}\sup_{\partial
B_{\nu r}}|v|
$$
\end{lemma}

Proof. Dilations show that it suffices to concentrate on
$r=1/\nu$. In that case the result follows from Theorem 5.5.8 of
\cite{Kr85} which states that
$$
|D^{2}v(x)-D^{2}v(y)|\leq N|x-y|^{\beta}\sup_{B_{1}}|v|
$$
as long as $x,y\in B_{1/2}$. The lemma is proved.

Introduce $\bL_{\delta}$ as the collection of operators
$Lu=a^{ij}D_{ij}u$ with $a=(a^{ij})$ being
measurable and $\cS_{\delta}$-valued.

\begin{lemma}
                                          \label{lemma 12.16.2}

Let $r\in(0,\infty)$ and let
 $w\in C^{2}(\bar{B}_{r})$ be a  function such that
$w=0$ on $\partial B_{r}$. 
Then there are   constants $\gamma\in(0,1]$ and $N$,
depending only on $\delta$ and $d$, such that
for any $L
\in\bL_{\delta}$  we have
$$
 \dashint_{B_{r}}|D^{2}w|^{\gamma}
\,dx \leq N\big(\dashint_{B_{r}}|Lw|^{d}\,dx\big)^{\gamma/d}.
$$
\end{lemma}

For $r=1$ the result is proved in
\cite{FH} on the basis of an
elliptic counterpart of Theorem \ref{theorem 12.21.1}.
 In the general case it suffices to use dilations.

\begin{lemma}
                                            \label{lemma 12.16.3}

Let   $r\in(0,\infty)$, $\nu\geq2$, and let
 $u\in W^{2}_{d}(B_{\nu r})$ be a   solution of
\eqref{12.16.3} in $B_{\nu r}$. Then
$$
\dashint_{B_{r}}
\dashint_{B_{r}}|D^{2}u(x)-D^{2}u(y)|^{\gamma}
\,dxdy 
$$
\begin{equation}
                                               \label{12.16.4}
 \leq  N\nu^{d}\big(\dashint_{B_{\nu
r}}\bar{f}^{d}\,dx\big)^{\gamma/d}
+N\nu^{ -\gamma\beta}  \big(\dashint_
{B_{\nu r}}
|D^{2}u|^{d}\,dx\big)^{\gamma/d},
\end{equation}
where $N=N(\delta,d)$.

\end{lemma}

Proof. Observe that it suffices to prove the lemma
for $u\in C^{\infty}_{b}(\bar{B}_{\nu r})$. Indeed,   if
we denote by $u^{n}$ a sequence converging to $u$
in $W^{2}_{d}(B_{\nu r})$, then
$$
\sup_{\omega}\big[a^{ij}(\omega,x)D_{ij}u^{n}(x)
+f^{n}(\omega,x)\big]=0
$$
in $B_{\nu r}$, where $f^{n}(\omega,x)=
a^{ij}(\omega,x)D_{ij}(u-u^{n})+f(\omega,x)$, so that
$$
\|\sup_{\omega}|f^{n}(\omega,\cdot)|\|_{\cL_{d}(B_{\nu r})}
\leq\|\bar{f}\|_{\cL_{d}(B_{\nu r})}+N
\|u-u^{n}\|_{W^{2}_{d}(B_{\nu r})},
$$
where the last term tends to zero as $n\to\infty$.

Next, let $v$ be a solution of \eqref{12.16.3} in $B_{\nu r}$
with $f^{\alpha}\equiv0$ and the boundary condition
$\hat{u}(x):=u(x)-x^{i}(u_{x^{i}})_{B_{\nu r}}-(u)_{B_{\nu r}}$ on
$\partial B_{\nu r}$. By classical results (see, for instance,
\cite{CC}, \cite{Ev}, \cite{GT}, \cite{Kr85}) such
a solution of class $C^{2}(\bar{B}_{\nu r})$ exists.
Then by Lemmas \ref{lemma 12.6.1} and~\ref{lemma 12.16.01}
\begin{equation}
                                                 \label{12.17.1}
 \dashint_{B_{r}}\dashint_{B_{r}}|D^{2}v(x)-D^{2}v(y)|^{\gamma}
\,dxdy \leq N\nu^{ -\gamma\beta}  \big(\dashint_
{B_{\nu r}}
|D^{2}u|^{d}\,dx\big)^{\gamma/d}.
\end{equation}

Furthermore, for $\hat{w}:=\hat{u}-v$ we have that in $B_{\nu r}$
$$
0=\sup_{\omega}\big[a^{ij}(\omega)D_{ij}\hat{w}(x)
+a^{ij}(\omega)D_{ij}v(x)+f(\omega,x)\big]
$$
$$
\leq \sup_{\omega}\big[a^{ij}(\omega)D_{ij}\hat{w}(x) +f(\omega,x)\big]
$$
and
$$
0=\sup_{\omega}\big[
 a^{ij}(\omega)D_{ij}\hat{w}(x)
+a^{ij}(\omega)D_{ij}v(x)+f(\omega,x)\big]
$$
$$
\geq \inf_{\omega}\big[a^{ij}(\omega)D_{ij}\hat{w}(x) +f(\omega,x)\big].
$$
It follows that there exists an operator $L\in\bL_{\delta}$ and a
function
$g$ such that $L\hat{w}+g=0$ in $B_{\nu r}$ and $|g|\leq\bar{f}$.
Therefore, by Lemma \ref{lemma 12.16.2} and by the fact that
$D^{2}\hat{w}=D^{2}w$, where $w=u-v$, we get that
$$
 \dashint_{B_{r}}|D^{2}w|^{\gamma}
\,dx \leq \nu^{d}\dashint_{B_{\nu r}}|D^{2}w|^{\gamma}
\,dx \leq
 N\nu^{d}\big(\dashint_{B_{\nu r}}\bar{f}^{d}\,dx\big)^{\gamma/d},
$$
$$
\dashint_{B_{r}}\dashint_{B_{r}}|D^{2}w(x)-D^{2}w(y)|^{\gamma}
\,dxdy \leq
 N\nu^{d}\big(\dashint_{B_{\nu r}}\bar{f}^{d}\,dx\big)^{\gamma/d}.
$$

By combining this with \eqref{12.17.1} we get
\eqref{12.16.4} and the lemma is proved.

Lemma \ref{lemma 12.16.3} plays the main role in the proof
of Theorem \ref{theorem 12.23.3}, a particular case of which
is the following theorem about apriori estimates and the
solvability of
Bellman's equations with constant coefficients.
\begin{theorem}
                                \label{theorem 12.19.1}
Let $p>d$.
(i) Let $u\in W^{2}_{p}$
satisfy \eqref{12.16.3} in $\bR^{d}$. Then there is a constant
$N=N(\delta,p,d)$ such that
\begin{equation}
                                           \label{12.19.5}
\|D^{2}u\|_{\cL_{p}}\leq N\|\bar{f}\|_{\cL_{p}}.
\end{equation}

(ii) There exists a $\lambda_{0}=\lambda_{0}
(\delta,p,K,d)\geq0$ such that for any $\lambda
\geq\lambda_{0}$ and
$u\in W^{2}_{p}$ satisfying
$$
\sup_{\omega\in \Omega}\big[a^{ij}(\omega )D_{ij}u(x)
+b^{i}(\omega,x)D_{i}u(x)
$$
\begin{equation}
                                           \label{12.19.2}
-(c+\lambda)(\omega,x)u(x)
+f(\omega,x)\big]=0
\end{equation}
in $\bR^{d}$ we have
\begin{equation}
                                           \label{12.19.3}
\lambda\|u\|_{\cL_{p}}+\|D^{2}u\|_{\cL_{p}}
\leq N\|\bar{f}\|_{\cL_{p}},
\end{equation}
where $N=N(\delta,p, d)$.
Moreover,  for any $\lambda
>\lambda_{0}$ there is a unique solution of \eqref{12.19.2}
in $W^{2}_{p}$. Finally, if $K=0$,
one can take $\lambda_{0}=0$.

\end{theorem}

Proof. (i) We use notation \eqref{12.19.4} for the filtration of
dyadic cubes and $\mu(dx)=dx$. We also note that one can apply
 Lemma \ref{lemma 12.16.3} after shifting the origin. Then we
easily  get that
 for any $\nu\geq2$
$$
(D^{2}u)_{\gamma}^{\sharp}
\leq N\nu^{d/\gamma}\bM^{1/d}(\bar{f}^{d})+
N\nu^{- \beta}\bM^{1/d}(|D^{2}u|^{d})
$$
on $\bR^{d}$, where $\bM g$ is the classical maximal function of $g$
and $N=N(\delta,d)$.

It follows by Theorem \ref{theorem 12.16.1}
and the Hardy-Littlewood theorem on maximal functions
that (recall that $p>d$)
\begin{equation}
                                               \label{12.16.9}
\|D^{2}u\|_{\cL_{p}}\leq N\nu^{d/\gamma}\|\bar{f}\|_{\cL_{p}}
+N\nu^{- \beta}\|D^{2}u\|_{\cL_{p}},
\end{equation}
where $N=N(\delta,p,d)$. The arbitrariness  of $\nu\geq2$
now leads to \eqref{12.19.5}.

(ii) By standard arguments all assertions in (ii)
follow from its last one, to prove which
it suffices to prove apriori estimate \eqref{12.19.3}.
Observe that by \eqref{12.19.5} we have (recall that $b=c=0$)
$$
\|D^{2}u\|_{\cL_{p}}\leq N\|\bar{f}\|_{\cL_{p}}
+N\lambda\|u\|_{\cL_{p}},
$$
which shows that it suffices to prove that
\begin{equation}
                                           \label{12.19.6}
\lambda\|u\|_{\cL_{p}} 
\leq N\|\bar{f}\|_{\cL_{p}}.
\end{equation}

Of course, we may assume that $\lambda>0$ and,
as in the proof of Lemma \ref{lemma 12.16.3}, we may assume that $u\in
C^{\infty}_{0}$. Then we take an operator
$L\in\bL_{\delta}$ and a function $g$ such that $|g|\le\bar{f}$ and
$Lu-\lambda u+g=0$. By Theorem 3.5.15 and Lemma
3.5.5 of \cite{Kr85} we have
\begin{equation}
                                           \label{12.19.7}
\lambda\|u\|_{\cL_{p}} \leq N\|g\|_{\cL_{p}}
\leq N\|\bar{f}\|_{\cL_{p}},
\end{equation}
with $N=N(\delta,d,p)$ provided that $\lambda\geq N$. Dilations
show that the latter requirement can be reduced to $\lambda>0$.
Actually, Theorem 3.5.15   of \cite{Kr85}
is only proved there for $p=d$, but the proof is based on the results
of Section 3.2 of  \cite{Kr85}, which treats general $p\geq d$
(see Theorem 3.2.3 there). Therefore, we can repeat
  the corresponding
arguments in the proofs of Theorem 3.5.15 and Lemma
3.5.5 of \cite{Kr85}
almost word for word and get what we need.
A different way to obtain the first estimate in \eqref{12.19.7}
is to refer directly to a rather old result of \cite{Kr72}
(see Lemma 7 there), but this would require the reader
to know It\^o's formula from the theory
of It\^o integrals. This way even has an advantage because it shows that
$N$ in \eqref{12.19.6} is independent of $p$.
The theorem is proved.

\mysection{Proof of Theorem \protect\ref{theorem 12.23.3}}
                                                \label{section 12.23.4}

We start with a few observations regarding
equations with variable~$a^{ij}$.
\begin{lemma}
                                            \label{lemma 12.19.1}

Let $\kappa\in(1,\infty)$,  $r\in(0,\infty)$, $\nu\geq2$, and let
 $u\in W^{2}_{d} $ be a   solution of
\eqref{12.6.1} in $\bR^{d}$. 
Assume that $u=0$ outside $B_{R_{0}}(x_{0})$ for some $x_{0}
\in\bR^{d}$. 
Then
$$
\dashint_{B_{r}}
\dashint_{B_{r}}|D^{2}u(x)-D^{2}u(y)|^{\gamma}
\,dxdy \leq
N\nu^{d}\big(\dashint_{B_{\nu
r}}\bar{f}^{d}\,dx\big)^{\gamma/d}
$$
\begin{equation}
                                               \label{12.16.7}
+ 
N\nu^{d}\big(\dashint_{B_{\nu r}}|D^{2}u|^{\kappa d}\,dx
\big)^{\gamma/(\kappa d)}\theta^{(1-1/\kappa)\gamma/d} 
+N\nu^{ -\gamma\beta}  \big(\dashint_
{B_{\nu r}}
|D^{2}u|^{d}\,dx\big)^{\gamma/d},
\end{equation}
where $N=N(\delta,d)$.

\end{lemma}

Proof. We basically repeat the proof of Lemma 6.2.2 of
\cite{Kr08}. Fix a $\nu\geq 2$ and an $r\in(0,\infty)$
and introduce
$$
\bar{a}^{ij}(\omega)=(a^{ij})_{B_{R_{0}}}\quad\text{if}\quad
\nu r\geq R_{0},\quad
\bar{a}^{ij}(\omega)=(a^{ij})_{B_{\nu r}}\quad\text{if}\quad
\nu r< R_{0}.
$$
Observe that
$$
\sup_{\omega\in \Omega}\big[\bar{a}^{ij}(\omega)D_{ij}u
+\tilde{f}(\omega,x)\big]=0,
$$
where
$$
\tilde{f}(\omega,x) =
(a^{ij}-\bar{a}^{ij})(\omega,x)D_{ij}u
+f(\omega,x).
$$
Denote
$$
\tilde{a}(x)=\sup_{\omega\in \Omega}|a(\omega,x)-\bar{a}(\omega)|.
$$
By Lemma \ref{lemma 12.16.3}
$$
\dashint_{B_{r}}
\dashint_{B_{r}}|D^{2}u(x)-D^{2}u(y)|^{\gamma}
\,dxdy  \leq  N\nu^{d}
I^{\gamma/d}
$$
\begin{equation}
                                               \label{12.16.8}
+ N\nu^{d}\big(\dashint_{B_{\nu
r}}\bar{f}^{d}\,dx\big)^{\gamma/d}
+N\nu^{ -\gamma\beta}  \big(\dashint_
{B_{\nu r}}
|D^{2}u|^{d}\,dx\big)^{\gamma/d},
\end{equation}
where
$$
I= \dashint_{B_{\nu
r}}\tilde{a}^{d}|D^{2}u|^{d}\,dx 
= \dashint_{B_{\nu
r}}\tilde{a}^{d}|D^{2}u|^{d}I_{B_{R_{0}}}\,dx\leq
J^{1/\kappa}_{1}J^{1-1/\kappa}_{2} 
$$
with
$$
J_{1}=\dashint_{B_{\nu r}}|D^{2}u|^{\kappa d}\,dx,
\quad
J_{2}=\dashint_{B_{\nu r}}\tilde{a}^{\kappa d/(\kappa-1)}
I_{B_{R_{0}}}\,dx\leq N
\dashint_{B_{\nu r}}\tilde{a} 
I_{B_{R_{0}}}\,dx.
$$

If $\nu r\geq R_{0}$, then
$$
J_{2}\leq N(\nu r)^{-d}\int_{B_{R_{0}}}\tilde{a} \,dx
\leq N(\nu r)^{-d}R^{d}_{0}\dashint_{B_{R_{0}}}\tilde{a} \,dx
\leq N\theta.
$$
In case $\nu r<R_{0}$ we have
$$
J_{2}\leq N \dashint_{B_{\nu r}}\tilde{a} \,dx
\leq N\theta.
$$
Hence,
$$
I\leq N\big(\dashint_{B_{\nu r}}|D^{2}u|^{\kappa d}\,dx
\big)^{1/\kappa}\theta^{1-1/\kappa}
$$
and by combining this with  \eqref{12.16.8}
we come to \eqref{12.16.7}. The lemma is proved.

\begin{corollary}
                                           \label{corollary 12.16.01}
Under the assumptions of Lemma \ref{lemma 12.19.1}
let $p>\kappa d$. Then there is a constant
$N_{1}=N_{1}(\delta,p,\kappa,d)$ such that
$$
\|D^{2}u\|_{\cL_{p}}\leq N_{1}
\nu^{ d/\gamma}\|\bar{f}\|_{\cL_{p}}
+N_{1}(\nu^{d/\gamma}\theta^{(1-1/\kappa) /d}+\nu^{-\beta})
\|D^{2}u\|_{\cL_{p}}
$$

\end{corollary}

This is obtained in the same way as \eqref{12.16.9}.

\begin{corollary}
                               \label{corollary 12.16.2}
Let $p>d$ and
let  $u\in W^{2}_{p} $ be a   solution of
\eqref{12.6.1} in $\bR^{d}$. 
Assume that $u=0$ outside $B_{R_{0}}(x_{0})$ for some $x_{0}
\in\bR^{d}$. Then
there exist $\theta=\theta(p,d,\delta)\in(0,1]$
and $N=N(p,d,\delta)$ such that
if Assumption \ref{assumption 12.18.1} is satisfied with this 
$\theta$, then $\|D^{2}u\|_{\cL_{p}}\leq N\|\bar{f}\|_{\cL_{p}}$.

\end{corollary}

Indeed, it suffices to set $2\kappa=1+p/d$ and choose first  $\nu$ and
then
$\theta$ in such a way that
$$
N_{1}(\nu^{ d/\gamma}
\theta^{(1-1/\kappa) /d}+\nu^{-\beta})\leq 1/2.
$$

{\bf Proof of Theorem  \ref{theorem 12.23.3}}. We suppose that 
Assumption \ref{assumption 12.18.1} holds with
$\theta$ from Corollary \ref{corollary 12.16.2}. 
First assume that we are given a function $u\in W^{2}_{p}$,
which satisfies \eqref{12.6.1} in $\bR^{d}$. 

Take
a nonnegative $\zeta\in C^{\infty}_{0}$ which 
 has support in $B_{R_{0}}$
and is such that $\zeta^{p}$ integrates to one.
 For the parameter $x_{0}\in\bR^{d}$
define
$$
u_{x_{0}}(x)=u(x)\zeta(x-x_{0})
$$
and observe that
$$
\sup_{\omega\in \Omega}\big[a^{ij}(\omega,x)D_{ij}u_{x_{0}}(x)
+f_{x_{0}}(\omega,x)\big]=0,
$$
where
$$
f_{x_{0}}(\omega,x)=f (\omega,x)\zeta(x-x_{0})
-u(x)a^{ij}(\omega,x)D_{ij}\zeta(x-x_{0})
$$
$$
-
2a^{ij}(\omega,x)(D_{i}u(x))D_{j}\zeta(x-x_{0}).
$$

By Corollary \ref{corollary 12.16.2}
$$
\|\zeta(\cdot-x_{0})|D^{2}u|\|_{\cL_{p}}^{p}
\leq N\|\zeta(\cdot-x_{0})\bar{f}\|_{\cL_{p}}^{p}
$$
$$
+\||D\zeta(\cdot-x_{0})|\,|Du|\|_{\cL_{p}}^{p}
+\||D^{2}\zeta(\cdot-x_{0})|u\|_{\cL_{p}}^{p}.
$$
Upon integrating through this estimate we get
\begin{equation}
                                              \label{12.16.10}
\| D^{2}u \|_{\cL_{p}}^{p}
\leq N_{1}\| \bar{f}\|_{\cL_{p}}^{p}
+N_{2}(\|   Du \|_{\cL_{p}}^{p}
+\| u\|_{\cL_{p}}^{p}),
\end{equation}
where $N_{1}=N_{1}(p,d,\delta)$ and $N_{2}=N_{2}
(p,d,\delta, R_{0})$.

If $u\in W^{2}_{p}$ is a solution of \eqref{12.23.5},
then by absorbing the first- and zeroth-order terms
into $f$ we see that
$$
\| D^{2}u \|_{\cL_{p}} 
\leq N_{1}(\| \bar{f}\|_{\cL_{p}} +\lambda\| u\|_{\cL_{p}})
+N_{2}(\|   Du \|_{\cL_{p}} 
+\| u\|_{\cL_{p}} )
$$
and if $\lambda\geq\lambda_{0} (p,d,\delta, R_{0})$, then
multiplicative inequalities yield
$$
\| D^{2}u \|_{\cL_{p}} 
\leq N_{1} (\| \bar{f}\|_{\cL_{p}} +\lambda\| u\|_{\cL_{p}}),
$$
where $N_{1}$ still depends only on $ p,d,\delta $.
Applying the results of \cite{Kr85} or \cite{Kr72}
as in the proof of Theorem \ref{theorem 12.19.1}, we obtain that
$$
\lambda\| u\|_{\cL_{p}}+
\| D^{2}u \|_{\cL_{p}} 
\leq N_{1}  \| \bar{f}\|_{\cL_{p}}.
$$
After that it suffices to repeat the corresponding argument
from the proof of Theorem \ref{theorem 12.19.1}.

\mysection{Comments on parabolic equations}

                                                \label{section 12.23.1}

Denote by $\frL$ the set of operators $L$ of the form
$$
L=\partial_{t}+a^{ij}(t,x)D_{ij}+b^{i}(t,x)D_{i}-c(t,x),\quad
\partial_{t}=\frac{\partial}{\partial t},
$$
where $a(t,x)=(a^{ij}(t,x))$ is an $\cS_{\delta}$-valued,
$b(t,x)=(b^{i}(t,x))$ is an $\bR^{d}$-valued,
and $c(t,x)$ is a real-valued measurable functions
defined on $\bR^{d+1}=\{(t,x):t\in\bR,x\in\bR^{d}\}$
satisfying 
$$
|b |+c \leq K,\quad c \geq0.
$$
Let $\frL_{0}$ be a subset of $\frL$ consisting of operators
with infinitely differentiable coefficients.

For $\rho,r>0$ introduce 
$$
C_{\rho,r}=(0,\rho)\times B_{r},\quad
\partial'C_{\rho,r}=\big([0,\rho]\times\partial B_{r}\big)
\cup\big(\{\rho\}\times B_{r}\big),
$$
$$
C_{\rho,r}(t,x)=(t,x)+C_{\rho,r},\quad
\partial'C_{\rho,r}(t,x)=(t,x)+\partial'C_{\rho,r}.
$$

\begin{theorem}
                                           \label{theorem 12.21.1}
Let $u\in W^{1,2}_{d+1}(C_{2,1})$ and assume that
$u\geq0$ on $\partial'C_{2,1}$ and there exists
an operator $L\in\frL$ such that $Lu\leq0$ in $C_{2,1}$.
Then there exist  constants $\gamma=\gamma(\delta,d,K)\in(0,1)$
and $N =N (\delta,d,K)$
such that for any $\lambda>0$
\begin{equation}
                                                    \label{12.21.2}
|C_{1,1}(1,0)\cap\{-Lu\geq\lambda\}| 
\leq N \lambda^{-\gamma}u^{\gamma}(0,0).
\end{equation}

\end{theorem}

Here is a consequence of this theorem, which can be
used in constructing the theory of parabolic Bellman's equations
along the lines in Sections \ref{section 12.23.3}
 and \ref{section 12.23.4}.

\begin{corollary}
                                             \label{corollary 12.23.1}

Let
 $w\in C^{1,2}(\bar{C}_{2,1})$ be a  function such that
$w=0$ on $\partial'C_{2,1}$. 
Then there are   constants $\gamma\in(0,1]$ and $N$,
depending only on $\delta$, $K$, and $d$, such that
for any $L
\in\frL$  we have
$$
 \int_{C_{1,1}}|D^{2}w|^{\gamma}
\,dxdt \leq N\big(\int_{C_{2,1}}|Lw|^{d+1}\,dxdt\big)^{\gamma/(d+1)}.
$$
\end{corollary}
This corollary is deduced from Theorem \ref{theorem 12.21.1}
in the same way as the theorem in \cite{FH} is deduced from
estimate (2.1) of \cite{FH}, the only difference being
that instead of the elliptic Alexandrov estimate one uses
the parabolic one.

To prove Theorem \ref{theorem 12.21.1} we need an auxiliary
construction, but first we observe that
replacing $u$ with $u/\lambda$ reduces the general case
to the one with $\lambda=1$.

For $q\in[0,1]$ denote by $\frU_{q}$ the set
of functions $u$, which are bounded and continuous
along with $\partial_{t}u$, $Du$, $D^{2}u$ in
$\bar{C}_{2,1}$ and such that

(i) $u=0$ on $\partial'C_{2,1}$;

(ii) there exists an operator 
$L\in \frL_{0}$ such that $Lu\leq0$ in $C_{2,1}$ and
$$
|C_{1,1}(1,0)\cap\{Lu\leq-1\}|\geq q|C_{1,1}|.
$$
Finally introduce
$$
m(t,x,q)=\inf\{u(t,x):u\in\frU_{q}\},\quad (t,x)\in \bar{C}_{2,1}.
$$

\begin{remark}
                                          \label{remark 12.21.1}
Denote by $m'(t,x,q)$ the function
called $m(t,x,q)$ in Section 4.1 of \cite{Kr85}. Then
obviously  $m'(t,x,q)=m(1+t,x,q)$. In particular,
by Lemmas 4.1.3 and 4.1.4 of \cite{Kr85}
for any $\kappa\in(0,1)$ there exist  $q_{0}\in(0,1)$ (close to 1)
and $m_{0}>0$, depending only on
$\delta,d,K$, and $\kappa$ such that
$$
m(t,x,q)\geq m_{0}
$$
for $q\in[q _{0},1]$ and $(t,x)\in\bar{C}_{\kappa^{2},\kappa}$.

\end{remark}

\begin{remark}
                                          \label{remark 12.21.2}
By Lemma 4.1.1 of \cite{Kr85} we know that
for $u$ from Theorem \ref{theorem 12.21.1} it holds that
$$
u(0,0)\geq m(0,0,q)  
$$
if $q$ satisfies
$$
q|C_{1,1}|=|C_{1,1}(1,0)\cap\{Lu\leq-1\}|.
$$
Therefore, to prove the theorem, we only need to prove that
\begin{equation}
                                              \label{12.21.3}
 N m(0,0,q)  \geq q^{1/\gamma}
\end{equation}
for $N$ and $\gamma$ depending only on $d,\delta,K$.
\end{remark}

We also need to use  
 Lemma 4.1.6 of \cite{Kr85}. Before stating
it we introduce some notation. Let $\Gamma$ be a measurable subset
of $C_{1,1}(1,0)$ and let $q,\eta,\zeta\in(0,1)$ be some numbers.
Denote by $\frB=\frB(\Gamma,q)$
 the collection of $Q=C_{\rho^{2},\rho}(t_{0},x_{0})$
such that $Q \subset C_{1,1}(1,0)$ and
$$
|Q\cap\Gamma|\geq q|C_{1,1}|.
$$
If $Q=C_{\rho^{2},\rho}(t_{0},x_{0})\in\frB$ we set
$$
Q'=(t_{0},x_{0})-C_{\eta^{-1}\rho^{2},\rho},\quad
Q''=(t_{0}-\eta^{-1}\rho^{2},x_{0})+C_{\eta^{-1}\rho^{2}
\zeta^{2},\rho\zeta}.
$$
Imagine that the $t$-axis is pointed up vertically. Then 
 $Q'$ is immediately adjacent to $Q$ from below, the two cylinders have a common base, and along the $t$-axis
$Q'$ is $\eta^{-1}$ times longer than $Q$. It is quite possible that
part of $Q'$ comes out of $C_{1,1}(1,0)$ or, for that matter, 
out of $C_{2,1}$. The cylinder $Q''$ is obtained from 
$ Q'$
by parabolic compression centered at the center of 
the lower base of $Q'$, the compression coefficient being $\zeta^{-1}$.
Finally, denote
$$
\Gamma''=\bigcup_{Q\in\frB}Q''.
$$
Here is   Lemma 4.1.6 of \cite{Kr85}.
\begin{lemma}
                                           \label{lemma 12.21.3}
If $|\Gamma|\leq q|C_{1,1}|$, then
$$
|\Gamma''|\geq \big(1-(1-q)3^{-d-1}\big)^{-1}(1+\eta)^{-1}
\zeta^{d+2}|\Gamma|.
$$

\end{lemma}

{\bf Proof of Theorem \ref{theorem 12.21.1}}.
We are going to slightly modify the
proof of Theorem 4.1.2 of \cite{Kr85} while concentrating on
\eqref{12.21.3}. Fix some $\eta,\zeta\in(0,1)$ to be specified later
and such that
\begin{equation}
                                                    \label{12.21.7}
\zeta^{2}\leq1-2\eta.
\end{equation}
Next, fix a $\kappa\in(0,1)$ such that $\kappa^{2}\geq1/2$ and set
$$
\mu(q)=\inf\{m(0,x,q):|x|\leq\kappa\}.
$$
By Remark \ref{remark 12.21.1} there is a $q_{0}\in(0,1)$
and $m_{0}>0$ such that
$$
\mu(q_{0})\geq m_{0}.
$$

Next, we take some $0<q'<q''<1$
and try to relate $\mu(q')$ to $\mu(q'')$.
To this end take a $u\in \frU_{q'}$ and let 
$$
\Gamma= C_{1,1}(1,0)\cap\{Lu\leq-1\}
$$
where $L$ is the operator associated with $u$. 
From chosen $\Gamma,q_{0},\eta$, and $\zeta$ we construct the
set $\Gamma''$ as before Lemma \ref{lemma 12.21.3} by taking there
$q_{0}$ in place of $q$ and consider two cases:

(i) $|\Gamma''\setminus C_{1,1}(1,0)|\leq (q''-q')|C_{1,1}|$,

(ii) $|\Gamma''\setminus C_{1,1}(1,0)|> (q''-q')|C_{1,1}|$.

{\em Case (i)\/}.  Denote by $\tilde{u}$ and $\tilde{v}$
the $W^{1,2}_{d+1}(C_{2,1})$-solutions of
$$
L\tilde{u}=-I_{\Gamma},\quad L\tilde{v}=-I_{\Gamma_{0}}
$$
vanishing on $\partial'C_{2,1}$, where $\Gamma_{0}=\Gamma''\cap
C_{1,1}(1,0)$. 
Since the coeficients of $L$ are infinitely differentiable,
such solutions exist.
There are two possibilities:  either

(a) $|\Gamma|\geq q_{0}|C_{1,1}|$,

\noindent or

(b) $|\Gamma|< q_{0}|C_{1,1}|$.

Under condition (a) by definition
\begin{equation}
                                                    \label{12.24.1}
u(0,x)\geq \mu(q_{0}),\quad|x|\leq\kappa.
\end{equation}

In case (b) by definition and  Lemma \ref{lemma 12.21.3}
$$
q'|C_{1,1}|\leq|\Gamma|\leq \big(1-(1-q_{0})3^{-d-1}\big)(1+\eta) 
\zeta^{-d-2}|\Gamma''|.
$$
Moreover, by assumption
$$
|\Gamma''|=|\Gamma''\setminus C_{1,1}(1,0)|+|\Gamma_{0}|
\leq (q''-q')|C_{1,1}|+|\Gamma_{0}|.
$$
It follows that
$$
|\Gamma_{0}|\geq q''|C_{1,1}|,
$$
if  
\begin{equation}
                                               \label{1.5.1}
(1+\xi) q'\geq 2q'',
\end{equation}
where
$$
\xi:=(1-(1-q_{0})3^{-d-1}\big)^{-1}(1+\eta)^{-1} 
\zeta^{ d+2} .
$$
Obviously there exist $\eta=\eta(q_{0})\in(0,1)$
and $\zeta=\zeta(q_{0})\in(0,1)$ such that
\eqref{12.21.7} is satisfied and $\xi=\xi(q_{0})>1$,
so that \eqref{1.5.1} holds for some $q'<q''$.
Since $q_{0}$  depends only on $\delta,K,d$, and $\kappa$,
so do $\eta$, $\zeta$,  and $\xi$.
We fix such $\eta$ and $\zeta$ from this moment on.
Then by definition
$$
\tilde{v}(0,x)\geq\mu(q'') ,\quad|x|\leq\kappa.
$$

We thus have estimated $\tilde{v}$ from below in case (i), (b)
 for $q'<q''$ satisfying \eqref{1.5.1}.
 By the maximum
principle $u\geq\tilde{u}$ and to estimate $u$ from below
it suffices to estimate $\tilde{u}$ from below
in terms of $\tilde{v}$. This will be done
by use of Lemma 4.1.5 of \cite{Kr85}.

If $(t_{0},x_{0})\in\Gamma_{0}$, then there exists
a cylinder $Q\in\frB$ such that $(t_{0},x_{0})\in Q''$.
Define $(t_{1},x_{1})$ and $\tau$, $\rho$ from the equation
$$
C_{\tau,\rho}(t_{1},x_{1})=Q_{1}:=(Q'\cup Q)\cap C_{2,1} ,
$$
so that
$$
Q=C_{ \rho^{2},\rho}(\tau+t_{1}-\rho^{2},x_{1}).
$$
Furthermore, owing to \eqref{12.21.7},
 the distance from $(t_{0},x_{0})$ to the bottom
of $Q$ is bigger than
$$
\eta^{-1}\rho^{2}-\eta^{-1}\rho^{2}\zeta^{2}\geq 2\rho^{2}.
$$
In particular,
\begin{equation}
                                               \label{12.21.8}
1<t_{0}\leq (\tau+t_{1}-\rho^{2})-2\rho^{2}.
\end{equation}

Next  let $\bar{u}$ and $\bar{v}$ be the 
$W^{1,2}_{d+1}(Q_{1})$-solutions of
$$
L\bar{u}=-N_{0}I_{\Gamma},\quad L\bar{v}=-I_{\Gamma_{0}},
$$
vanishing on $\partial'Q_{1}$, where the constant $N_{0}$
will be specified later in such a way that
\begin{equation}
                                               \label{12.21.5}
\bar{u}(t_{0},x_{0})\geq \bar{v}(t_{0},x_{0}).
\end{equation}
If we can do this, then by Lemma 4.1.5 of \cite{Kr85}
we have $N_{0}\tilde{u}\geq\tilde{v}$ on $\bar{C}_{2,1}$ and
  for $q'<q''$ satisfying \eqref{1.5.1}
\begin{equation}
                                               \label{12.21.6}
u(0,x)\geq\tilde{u}(0,x)\geq N_{0}^{-1}\tilde{v}(0,x)
\geq N_{0}^{-1}\mu(q''),\quad|x|\leq\kappa.
\end{equation}

To prove \eqref{12.21.5}, observe that by the maximum 
principle $\bar{v}(t,x)\leq t_{1}+\tau-t$, so that
\begin{equation}
                                               \label{12.22.1}
\bar{v}(t_{0},x_{0})\leq(1+\eta^{-1})\rho^{2}.
\end{equation}
  On the other hand, by the choice of $Q$ we have
$q_{0}|Q|\leq|\Gamma\cap Q|$. This inequality is preserved
under the parabolic dilation $x\to\rho^{-1}x$, $t\to\rho^{-2}t$
which transforms $\bar{u}$ into a function $\hat{u}$,
which satisfies $\hat{L}\hat{u}=-N_{0}\rho^{2}I_{\hat{\Gamma}}$
with an $\hat{L}\in\frL_{0}$ and $\hat{\Gamma}$ being the
image of $\Gamma$. Observe that the hyperplane
$t=\tau+t_{1}-2\rho^{2}$ passes at a distance $\rho^{2}$
from $Q$ and it intersects $Q_{1}$ above $t=t_{0}$.
Hence, by definition
$$
\bar{u}(\tau+t_{1}-2\rho^{2},x) 
\geq N_{0}\rho^{2}\mu(q_{0}),\quad|x-x_{1}|\leq\kappa\rho.
$$
Since $|x_{0}-x_{1}|\leq (1-\zeta)\rho$ and
 the distance between the hyperplanes $t=\tau+t_{1}-2\rho^{2}$
and $t=t_{0}$ is bigger than $\rho^{2}$
and less than $\eta^{-1}\rho^{2}$, it follows 
 by Lemma 4.1.3 of \cite{Kr85} that
$$
\bar{u}(t_{0},x_{0}) 
\geq \alpha N_{0}\rho^{2}\mu(q_{0}), 
$$
where $\alpha>0$ depends only on $d,\delta,K$, and $\kappa$.
By taking
$$
N_{0}=\alpha^{-1}(1+\eta^{-1})\mu^{-1}(q_{0})
$$
and recalling \eqref{12.22.1}
we come to \eqref{12.21.5} and thus \eqref{12.21.6}
is established in case (i), (b)
 for $q'<q''$ satisfying \eqref{1.5.1}, so that
generally  in case (i) (recall \eqref{12.24.1})
 for those $q',q''$
$$
u(0,x)\geq\min(\mu(q_{0}), N^{-1}_{0}\mu(q'')),\quad|x|\leq\kappa.
$$
The arbitrariness in the choice of $u$ implies that
$$
\mu(q')\geq\min(\mu(q_{0}), N^{-1}_{0}\mu(q'')),
$$
which after introducing
$$
\hat{\mu}(q)=\min(\mu(q_{0}),\mu(q))
$$
yields
\begin{equation}
                                                    \label{12.24.2}
\hat{\mu}(q')\geq  N^{-1}_{0}\hat{\mu}(q''))
 \quad\text{if}\quad (1+\xi)q'\geq2q''.
\end{equation}

{\em Case (ii)\/}. First we claim that for some $(t,x)\in\Gamma''$
it holds that $t<q'-q''+1$. Indeed, otherwise
$\Gamma''\in C_{q''-q',1}(1,0)$ and $|\Gamma''|\leq
(q''-q')|C_{1,1}|$. It follows that there is a cylinder
$$
Q=C_{\rho^{2},\rho}(t_{0},x_{0})\in\frB
$$
 such that $Q'$
contains points in the half-space $t<q'-q''+1$. 
Since $q'<q''$, $q'-q''+1<1$ and $Q'$ is adjacent
to $Q\in C_{1,1}(1,0)$, this implies that the height of $Q'$
is at least $q''-q'$, that is,
\begin{equation}
                                            \label{12.22.7}
\rho^{2}\eta^{-1}\geq q''-q',\quad\rho^{2}\geq\eta(q''-q').
\end{equation}

Moreover, by definition  $|\Gamma\cap Q|\geq q_{0}|Q|$ 
and by using dilations and the maximum principle we see that
\begin{equation}
                                            \label{12.22.4}
u(t_{0}-\rho^{2},x)\geq \mu(q_{0})\rho^{2},\quad |x-x_{0}|\leq
\kappa\rho.
\end{equation}
If $t_{0}-\rho^{2}\geq1/4$, then \eqref{12.22.4}
 by Lemma 4.1.3 of \cite{Kr85} 
implies that
$$
u(0,x)\geq \mu(q_{0})\varepsilon\rho^{n} ,\quad |x|\leq
\kappa ,
$$
where $\varepsilon>0$ and $n\geq1$ depend 
only on $\delta,K,d$, and $\kappa$.
In particular (see \eqref{12.22.7}),
\begin{equation}
                                            \label{12.22.6}
u(0,x)\geq \mu(q_{0})\varepsilon\eta^{n}
(q''-q')^{n},\quad |x|\leq
\kappa .
\end{equation}

On the other hand, 
by Remark \ref{remark 12.21.1} for $t_{0}-\rho^{2}
\leq t\leq t_{0}-\rho^{2}+\kappa^{2}\rho^{2}$ and
$|x-x_{0}|\leq\kappa \rho$ it holds that
$$
u(t,x)\geq m_{0}\rho^{2}.
$$
In addition, if $t_{0}-\rho^{2}<1/4$, then
in the above inequality one can take $t
=\kappa^{2}/2$ since 
$t_{0}\geq1$, $\rho^{2}>3/4$, $\rho^{2}\leq1$
(also recall that $\kappa^{2}\geq1/2$) and 
$$
t_{0}-\rho^{2}<1/4\leq\kappa^{2}/2\leq t_{0}-\rho^{2}+\kappa^{2}/2
\leq t_{0}-\rho^{2}+\kappa^{2}\rho^{2}.
$$
Hence, $u(\kappa^{2}/2,x)\geq m_{0}\rho^{2}$ for $|x-x_{0}|\leq
\kappa\rho$ and as above
$$
u(0,x)\geq m_{0}\varepsilon\eta^{n}(q''-q')^{n},
\quad |x|\leq\kappa,
$$
with $\varepsilon$ and $n$ of the same kind as in
\eqref{12.22.6} or just the same if we choose
the minimum of $\varepsilon$'s and the maximum of $n$'s.
Again the arbitrariness of $u$ yields $\mu(q')\geq
m_{0}\varepsilon\eta^{n}(q''-q')^{n}$, which after
reducing $\varepsilon$ if necessary, so that
$\mu(q_{0})\geq m_{0}\varepsilon\eta^{n}$ leads to
$$
\hat{\mu}(q')\geq m_{0}\varepsilon\eta^{n}(q''-q')^{n}.
$$

As a result of considering the two cases (i) and (ii) we get that
there exist $\varepsilon_{0}\in(0,1)$ and $n_{0}\geq1 $
depending only on  $\delta,K,d$, and $\kappa$, such that
for any $0<q'<q''<1$ 
  such that $(1+\xi)q'\geq2q''$ we have
\begin{equation}
                                                 \label{12.23.2}
\hat{\mu}(q')\geq\varepsilon_{0}\min\big( 
(q''-q')^{n_{0}},\hat{\mu}(q'')\big).
\end{equation}
We also know that $\hat{\mu}(q)\geq m_{0}>0$ for $q\geq q_{0}$.

 We may certainly assume that $\varepsilon_{0}
\leq\bar{\varepsilon}:=2/(1+\xi)$ (recall that
$\xi>1$) and 
 we claim that for $q_{k}=
 \bar{\varepsilon}^{k} q_{0}$,
$k=0,1,2,...$, we have
\begin{equation}
                                                 \label{12.23.1}
\hat{\mu}(q_{k})\geq\varepsilon^{kn_{0}}_{0}\chi,\quad
\chi:=
\min\big(\mu(q_{0}),q_{0}^{n_{0}}(1-
 \bar{\varepsilon}^{n_{0}}
\big).
\end{equation}
To prove the claim we use induction. If $k=0$,
\eqref{12.23.1} is obvious. If it is true for a $k$,
then $q_{k}-q_{k+1}= \bar{\varepsilon}^{k}
q_{0}(1- \bar{\varepsilon})$
$$
(q_{k}-q_{k+1})^{n_{0}}= \bar{\varepsilon}^{kn_{0}}
q_{0}^{n_{0}}(1- \bar{\varepsilon})^{n_{0}}
\geq \varepsilon_{0}^{kn_{0}}\chi,
$$
so that by \eqref{12.23.2}  
and the fact that $(1+\xi)q_{k+1}=2q_{k}$
$$
\hat{\mu}(q_{k+1})\geq\varepsilon_{0}\min\big(\varepsilon_{0}^{kn_{0}}
\chi,\hat{\mu}(q_{k})\big)\geq\varepsilon_{0}\varepsilon_{0}^{kn_{0}}
\chi\geq \varepsilon_{0}^{(k+1)n_{0}}
\chi.
$$

This proves \eqref{12.23.1} and shows that,
  if we define $r>1$ so that $\varepsilon_{0}^{n_{0}}
=\bar{\varepsilon}^{r}$, then
$
\hat{\mu}(q_{k})\geq Nq_{k}^{r}$ with 
$r,N>0$ depending only on
 $\delta,K,d$, and $\kappa$. By observing that $\hat{\mu}$ 
is an increasing function we obtain that $\hat{\mu}(q)
\geq Nq^{r}
$, $\mu(q)\geq Nq^{r}$ for $q\leq1$. Finally, since
 $m(0,0,q)\geq\mu(q)$ if in the construction of $\mu$, we 
take any $\kappa>0$, say $\kappa^{2}=1/2$, we come to
\eqref{12.21.3} with $\gamma=1/r$,
which, as it is explained in
Remark \ref{remark 12.21.2}, proves the theorem.

\mysection{Appendix}
                                           \label{section 12.24.3}

We will be working in the setting of Chapter 3 of \cite{Kr08}
using notation different from the previous sections of 
the present article. Thus, $(\Omega,\cF,\mu)$ is a measurable space
with $\sigma$-finite measure $\mu$ satisfying
$\mu(\Omega)=\infty$. For $\Gamma\in\cF$ we use the notation
$|\Gamma|=\mu(\Gamma)$ and
$$
f_{\Gamma}=\dashint_{\Gamma}f\,\mu(dx)=\frac{1}{|\Gamma|}
\int_{\Gamma}f\,\mu(dx).
$$

Next we take a filtration $\{\bC_{n}:n\in\bZ\}$ of partitions
of $\Omega$ as in Section 3.1 of \cite{Kr08} and 
recall that for any $n\in\bZ$ and $C\in \bC_{n}$ there exists
a  unique ``parent'' $C'\in\bC_{n-1}$ such that $C\subset
C'$. It is assumed that whenever $C$ and $C'$ are related in the above
described manner we have $|C'|\leq N_{0}|C|$, where $N_{0}$
is a constant independent of $n$, $C$, and $C'$.

For functions $g$, for which it makes sense,
denote
$$
g_{|n}(x)=\dashint_{C_{n}(x)}g(y)\,\mu(dy),
$$
where $C_{n}(x)$ is the element of the family $\bC_{n}$
containing $x$. Also
$$
\cM g(x):=\sup_{n<\infty}|g|_{|n}(x).
$$

The most relevant filtration of partitions in this paper
is the dyadic cube filtration of partitions of $\bR^{d}$
with Lebesgue measure when
$$
\bC_{n}=\{C_{n}(i_{1},...,i_{d}),i_{1},...,i_{d}\in\bZ\},
$$
$$
C_{n}(i_{1},...,i_{d})=[i_{1}2^{-n},(i_{1}+1)2^{-n})\times...\times
[i_{d}2^{-n},(i_{d}+1)2^{-n}).
$$

In the remaining part of the section we consider two functions
$u,v\in\cL_{1}(\Omega)$ and a nonnegative measurable function $g$
on $\Omega$.   Below by $I_{\cM v(x)>\alpha \lambda}$
we mean the indicator function of the set
$\{x:\cM v(x)>\alpha \lambda\}$. 

The most relevant case of Lemma \ref{lemma 12.16.1}
for the purposes of the present article is
when $u^{C}=|u|$. In the form it is stated
and for $\gamma=1$ the lemma was used in \cite{Kr09}
while treating linear elliptic  equations with rather rough coefficients.

\begin{lemma}
                                         \label{lemma 12.16.1}
Let $\gamma\in(0,1]$. Assume that
  $|u|\leq v$
and 
for any $n\in\bZ$ and $C\in\bC_{n}$ there exists a
measurable function  $u^{C}$ given on $C$ such that
$|u|\leq u^{C}\leq v$ on $C$ and
\begin{equation}
                                                      \label{6.29.3}
 \dashint_{C}
\dashint_{C}|u^{C}(x)-u^{C}(y)|^{\gamma}\,\mu(dx)\mu(dy)\leq
 \dashint_{C}g^{\gamma}(x)\,\mu(dx) .
\end{equation}
Then for any $\lambda>0$ we have
\begin{equation}
                                                      \label{6.29.1}
 |\{x:|u(x)|\geq\lambda\}|\leq \nu^{-1}\lambda^{-\gamma}
\int_{\Omega}g^{\gamma}(x)I_{\cM v(x)>\alpha \lambda}\,\mu(dx),
\end{equation}
where $\alpha=(2N_{0})^{-1}$ and $\nu= 1-2^{-\gamma}$.
\end{lemma}

Proof. Obviously we may assume that   $u\geq0$.
Fix a $\lambda>0$ and define 
$$
\tau(x)=\inf\{n\in\bZ:v_{|n}(x)>\alpha\lambda\}.
$$
We know that $\tau$ is a stopping time and if $\tau(x)<\infty$,
then 
$$
  v_{|n}(x)\leq \lambda/2,\quad\forall n\leq\tau(x).
$$
We also know that $v_{|n}\to v\geq u$ (a.e.) as $n\to\infty$.
It follows that (a.e.)
$$
\{x:u(x)\geq\lambda\}=\{x:u(x)\geq\lambda,\tau(x)<\infty\}
$$
$$
=\{x:u(x)\geq\lambda,  v_{|\tau}(x)\leq \lambda/2\}
=\bigcup_{n\in\bZ}\bigcup_{C\in \bC^{\tau}_{n}}A_{n}(C),
$$
where
$$
A_{n}(C):=\{x\in C:u(x)\geq\lambda, v_{|n}(x)\leq \lambda/2\},
$$
and $\bC^{\tau}_{n}$ is the family of disjoint elements
of $\bC_{n}$ such that
$$
\{x:\tau(x)=n\}=\bigcup_{C\in \bC^{\tau}_{n}}C.
$$

Next, for each $n\in\bZ$ and $C\in\bC_{n}$ on the set $A_{n}(C)$,
if it is not empty,
we have $v_{|n}=v_{C}$ and on $A_{n}(C)$
$$
u^{\gamma}-(v_{C})^{\gamma}\geq\lambda^{\gamma}(1-2^{-\gamma})
=\nu\lambda^{\gamma} .
$$
  We use this
and the inequality $|a-b|^{\gamma}\geq|a|^{\gamma}-|b|^{\gamma}$
and conclude that for $x\in A_{n}(C)$
$$
\dashint_{C}|u^{C}(x)-u^{C}(y)|^{\gamma}\,\mu(dy)
\geq(u^{C}(x))^{\gamma}-\dashint_{C}(u^{C}(y))^{\gamma}\,\mu(dy)
$$
$$
\geq u^{\gamma}(x)-\dashint_{C}v^{\gamma}(y)\,\mu(dy)
\geq u^{\gamma}(x)-(v_{C}(x))^{\gamma}
\geq\nu\lambda^{\gamma},
$$
 so that
by Chebyshev's inequality 
$$
|A_{n}(C)|\leq 
\nu^{-1}\lambda^{-\gamma}\int_{C}
\dashint_{C}|u^{C}(x)-u^{C}(y)|^{\gamma}\,\mu(dy)\mu(dx).
$$

It follows by assumption \eqref{6.29.3} that
$$
|A_{n}(C)|\leq 
\nu^{-1}\lambda^{-\gamma}\int_{C}g^{\gamma} \,\mu(dx),
$$
$$
|\{x:u(x)\geq\lambda\}|\leq\nu^{-1}\lambda^{-\gamma}
\sum_{n\in\bZ}\sum_{C\in \bC^{\tau}_{n}}\int_{C}g ^{\gamma}\,\mu(dx)
$$
$$
=\nu^{-1}\lambda^{-\gamma}\int_{\Omega}g^{\gamma}I_{\tau<\infty}\,\mu(dx).
$$
It only remains to observe that $\{\tau<\infty\}=\{\cM
v>\alpha\lambda\}$. The lemma is proved.

\begin{corollary}
                                           \label{corollary 12.16.1}
Under the assumptions of Lemma \ref{lemma 12.16.1}
for any $p>\gamma$ we have
$$
\int_{\Omega}|u|^{p}\,\mu(dx)\leq\beta
\big(\int_{\Omega}(\cM v)^{p}\,\mu(dx)\big)^{1-\gamma/p}
\big(\int_{\Omega}g^{p}\,\mu(dx)\big)^{ \gamma/p},
$$
where $\beta=\nu^{-1}(1-\gamma/p)^{-1}\alpha^{\gamma-p}$.
\end{corollary}

Indeed,
$$
\int_{\Omega}|u|^{p}\,\mu(dx)=
\int_{0}^{\infty}|\{x:|u(x)|\geq\lambda^{1/p}\}|\,d\lambda
$$
$$
\leq\nu^{-1}\int_{\Omega}g^{\gamma}\int_{0}^{\infty}
\lambda^{-\gamma/p}I_{\cM v>\alpha\lambda^{1/p}}\,d\lambda
\,\mu(dx)
=\beta\int_{\Omega}g^{\gamma}(\cM v)^{p-\gamma}\,\mu(dx)
$$
and it only remains to use H\"older's inequality.

The second statement of the following theorem
for $\gamma=1$ is, actually, the Fefferman-Stein theorem.

\begin{theorem}
                                         \label{theorem 12.16.1}
For $\gamma\in(0,1]$  define
\begin{equation}
                                           \label{12.19.4}
u_{\gamma}^{\sharp}(x)=\sup_{n}\sup_{C\in\bC_{n}:x\in C}
\big( \dashint_{C}
\dashint_{C}|u (z)-u (y)|^{\gamma}\,
\mu(dz)\mu(dy)\big)^{1/\gamma}.
\end{equation}
Then for $p>\gamma$
$$
\int_{\Omega}|u|^{p}\,\mu(dx)\leq\beta
\big(\int_{\Omega}(\cM u )^{p}\,\mu(dx)\big)^{1-\gamma/p}
\big(\int_{\Omega}(u_{\gamma}^{\sharp})^{p}\,\mu(dx)\big)^{
\gamma/p}.
$$
In particular, by Hardy-Littlewood theorem, for $p>1$
and $u\in\cL_{p}$ we have
$$
\|u\|_{\cL_{p}}\leq
N\|u_{\gamma}^{\sharp}\|_{\cL_{p}}  ,
$$
where $N=\beta^{1/\gamma}q^{(p-\gamma)/\gamma}$,
$q=p/(p-1)$.
\end{theorem}

This is a simple consequence of Corollary \ref{corollary 12.16.1}
since, obviously, for $n\in\bZ$, $  C\in\bC_{n}$,
and $v=u^{C}=|u|$
we have $|u|\leq u^{C}\leq v$, and for any $x\in C$
$$
 \dashint_{C}
\dashint_{C}|u^{C}(z)-u^{C}(y)|^{\gamma}\,
\mu(dz)\mu(dy) 
$$
$$
\leq
 \dashint_{C}
\dashint_{C}|u (z)-u (y)|^{\gamma}\,
\mu(dz)\mu(dy)
\leq(u_{\gamma}^{\sharp}(x))^{\gamma},
$$
so that
$$
 \dashint_{C}
\dashint_{C}|u^{C}(z)-u^{C}(y)|^{\gamma}\,
\mu(dz)\mu(dy) \leq\dashint_{C}(u_{\gamma}^{\sharp} )^{\gamma}
\,\mu(dx).
$$

\end{document}